\documentclass[12pt]{article}
\usepackage[all]{xy}
\usepackage{amsmath}
\usepackage{amsfonts}
\usepackage{amssymb}
\usepackage{amscd}
\usepackage{amsthm}
\usepackage{latexsym}
\usepackage{amsbsy}

\newtheorem{lem}{Lemma}[section]
\newtheorem{prop}{Proposition}[section]
\newtheorem{theorem}{Theorem}[section]
\newtheorem{corollary}{Corollary}[section]

\begin{document}
\date{}
\title
{ Hopf Algebra Equivariant Cyclic Homology and  Cyclic Homology of Crossed Product Algebras}
\author{R. Akbarpour,  M. Khalkhali \\Department of Mathematics 
, University of Western Ontario\\
\texttt{akbarpur@uwo.ca \; masoud@uwo.ca}}
\maketitle

\begin{abstract}
We introduce the cylindrical module $A \natural \mathcal{H}$, where $\mathcal{H}$ is a Hopf algebra 
and $A$ is a Hopf module algebra over $\mathcal{H}$. We show that there exists an isomorphism
between $\mathsf{C}_{\bullet}(A^{op} \rtimes \mathcal{H}^{cop})$ the cyclic module of the crossed product algebra $A^{op} \rtimes \mathcal{H}^{cop} $, and $\Delta(A \natural \mathcal{H}) $, the cyclic 
module related to the diagonal of $A \natural \mathcal{H}$. If $S$, the antipode of $\mathcal{H}$, is invertible it follows that $\mathsf{C}_{\bullet}(A \rtimes \mathcal{H}) \simeq \Delta(A^{op} \natural \mathcal{H}^{cop})$. 
When $S$ is invertible, we approximate $HC_{\bullet}(A \rtimes \mathcal{H})$ by a spectral
sequence and  give an interpretation of $ \mathsf{E}^0 , \mathsf{E}^1$ and $\mathsf{E}^2 $ terms of this 
spectral sequence.  
\end{abstract}
\textbf{Keywords.}  Cyclic homology, Hopf algebras  .

%**********************************Introduction 1.***************************

\section{Introduction.}
A celebrated problem in cyclic homology theory is to compute the cyclic homology of the crossed product algebra $A \rtimes G $,
 where the group $G$ acts on the algebra $A$ by automorphisms. This problem has been studied
by E. Getzler and D.S. Jones~\cite{gj93}, Feigin and Tsygan~\cite{fs86}, J. L. Brylinski~\cite{b87} and V. Nistor~\cite{n89}.
If $G$ is a discrete group, there is a spectral sequence, due to Feigin and Tsygan~\cite{fs86}, which converges to the cyclic homology
of the crossed product algebra. This result generalizes Burghelea's calculation of the cyclic homology of a 
group algebra~\cite{b85}.
 
In their paper, Getzler and Jones introduced the cylindrical module $A \natural G$ and 
they showed that there is an isomorphism
of cyclic modules
$\Delta(A \natural G) \cong \mathsf{C}_{\bullet}(A \rtimes G)$, where $\Delta$ denotes the diagonal of the cylindrical 
module and $\mathsf{C}_{\bullet}$ denotes the cyclic module of an algebra. Using their Eilenberg-Zilber 
theorem for cylindrical modules,
, they were then able to derive the Feigin-Tsygan spectral sequence to compute cyclic homology of $A \rtimes G $.
In this paper, we generalize their work to the case of Hopf
module algebras with the action of a Hopf algebra.

In Section 2, we review some  facts about paracyclic modules, cylindrical modules and Eilenberg-Zilber theorem for cylindrical modules~\cite{gj93}, and we introduce the terminology that
we need for the other parts of the paper. For some references to this section one can see~\cite{gj93,fs86,sw69,ll92}.
In Section 3, we introduce the cylindrical module $A \natural \mathcal{H}$ and we 
conclude that its diagonal $\Delta(A \natural \mathcal{H})$
is a cyclic module.
In Section 4, we show that the cyclic modules $\mathsf{C}_{\bullet}(A^{op} \rtimes \mathcal{H}^{cop}) $ 
and $\Delta (A \natural \mathcal{H})$ 
are isomorphic. When $S$, the antipode of $\mathcal{H}$, is invertible, we can then 
conclude that $\mathsf{C}_{\bullet}(A \rtimes \mathcal{H}) $ and
 $\Delta(A^{op} \natural \mathcal{H}^{cop})$ are isomorphic.
In Section 5, we construct the equivariant cyclic module $C_{\bullet}^{\mathcal{H}^{cop}}(A^{op})$ and 
when $S$ is invertible $C_{\bullet}^{\mathcal{H}}(A)$, dual to the cyclic module introduced in~\cite{rm00}, and 
we approximate the cyclic homology of 
$\mathsf{C}_{\bullet}(A^{op}  \rtimes \mathcal{H}^{cop})$ 
and $\mathsf{C}_{\bullet}(A  \rtimes \mathcal{H})$(when again $S$ is invertible), by a spectral sequence. We also compute
the $\mathsf{E}^0$, $\mathsf{E}^1$, and $\mathsf{E}^2$ terms of this spectral sequence.

In an earliear version of this paper~\cite{rm01} we were only able to prove our main result under the assumption $S^2=id_{\mathcal{H}}$
for the antipode of $\mathcal{H}$. We find it remarkable that a modification of our original formulas work for all Hopf algebras with invertible antipode.
It is known that such Hopf algebras form a very large class that includes quantum groups and quantum enveloping algebras.

%***************************************************************
\section{Preliminaries on paracyclic modules and Hopf algebras   
}

Let $[j],  j \in \mathbb{Z}^{+} $ be the set of all integers with the standard ordering and 
$t_j$ an automorphism of $[j]$ defined by $ t_j(i)=i+1,i \in \mathbb{Z} $. The paracyclic
category $\Lambda_{\infty}$  is a category whose
objects are $[j], j \ge 0$ and whose morphisms, $Hom_{\Lambda_{\infty}}([i],[j])$, are the
sets of all nondecreasing maps $f:[i] \rightarrow [j]$ with $t_j^{j+1}f=f \; t_i^{i+1}$,~\cite{fs86}.

One can identify the simplicial category $\Delta$ as a sub category of $\Lambda_{\infty}$ so that the objects 
of $\Delta$ are the same as $\Lambda_{\infty}$ and morphisms, $Hom_{\Delta}([i],[j])$, are those morphisms of 
$Hom_{\Lambda_{\infty}}([i],[j])$ that map $\{0,1, \dots ,i \}$ into $\{ 0,1, \dots , j \}$.\\
Every element $g \in Hom_{\Lambda_{\infty}}([i],[j]) \; \; , 0 \le r \le \infty ,$ has a unique
decomposition $g=f\; t_i^k $ where $ f \in Hom_{\Delta}([i],[j])$.

If $\mathcal{A}$ is any category, a paracyclic object in $\mathcal{A}$ 
is a contravariant functor from $\Lambda_{\infty}$ to $\mathcal{A}$. This definition is equivalent to giving 
a sequence of objects $A_0 , A_1 , \dots$ together with face operators $\partial_i : A_n \rightarrow A_{n-1}$ 
, degenerecy operators $\sigma_i : A_n \rightarrow A_{n+1} \;\;(i=0,1,\dots ,n),\;$and cyclic operators 
$\tau_n:A_n \rightarrow A_n$ where these operators
satisfy the simplicial
and $\Lambda_{\infty}$-identities,      
\begin{eqnarray} \label{eq:sim} \notag
& &\partial_i\partial_j = \partial_{j-1}\partial_i \hspace{15pt} i<j\\ \notag
& &\sigma_i\sigma_j = \sigma_{j+1}\sigma_i \hspace{15pt} i\le j \\ 
& &\partial_i\sigma_j =
\begin{cases}
\sigma_{j-1}\partial_i &\text{$i<j$}\\
\text{identity}      &\text{$i=j$ or $i=j+1$}\\
\sigma_j\partial_{i-1} &\text{$i>j+1$}.
\end{cases}
\end{eqnarray} 
\begin{eqnarray*} \label{eq:pa}
\partial_i\tau_n = \tau_{n-1}\partial_{i-1} \hspace{15pt} 1\le i \le n\; , \quad
\partial_0 \tau_n = \partial_n \notag \\
\sigma_i \tau_n =  \tau_{n+1} \sigma_{i-1} \hspace{15pt} 1\le i\le n\; , \quad
\sigma_0 \tau_n = \tau_{n+1}^2 \sigma_n 
\end{eqnarray*}
If in addition we have $\tau_n^{n+1}=id$, then we have a cyclic object in the sense of Connes~\cite{aC994}.

By a bi-paracyclic object in a category $\mathcal{A}$, we mean a paracyclic object in
the category of paracyclic objects on $\mathcal{A}$. So, giving a bi-paracyclic object
on  $\mathcal{A}$ is equivalent to giving a double sequence $A(p,q)$ of elements of $\mathcal{A}$ and operators
$\partial_{(p,q)} ,\sigma_{(p,q)},\tau_{(p,q)}$ and
$ \bar{\partial}_{(p,q)} ,\bar{\sigma}_{(p,q)},\bar{\tau}_{(p,q)}$ such that, for all $p \ge 0$,
\begin{eqnarray*}
B_p(q)=\{ A(p,q),\sigma_i^{p,q} ,\partial_i^{p,q} ,\tau_{p,q} \;\;\; :0 \le i \le q \}
\end{eqnarray*}
and for all $q \ge 0$,
\begin{eqnarray*}
\bar{B}_q(p)=\{ A(p,q),\bar{\sigma}_i^{p,q} ,\bar{\partial}_i^{p,q} ,\bar{\tau}_{p,q} \;\;\; :0 \le i \le p \}
\end{eqnarray*}
are paracyclic objects in $\mathcal{A}$ and every horizontal operator commutes with every vertical operators.

We say that a bi-paracyclic object is cylindrical~\cite{gj93} if for all $p,q \ge 0,$
\begin{equation} \label{eq:cy}
\bar{\tau}_{p,q}^{q+1} \;
\tau_{p,q}^{p+1}=\tau_{p,q}^{p+1} \;
\bar{\tau}_{p,q}^{q+1}=id_{p,q}.
\end{equation}
If $A$ is a bi-paracyclic object in a category $\mathcal{A}$, the paracyclic
object related to the diagonal of $A$ will be denoted by $\Delta A$. So the paracyclic operators on $\Delta A(n)=A(n,n)$ are $\bar{\partial}_i^{n,n-1}
\; \partial_i^{n,n} , \bar{\sigma}_i^{n,n+1} 
\sigma_i^{n,n},\bar{\tau}_{n,n}  \tau_{n,n}$, where $0 \le i \le n.$
When $A$ is a cylindrical object since the cyclic operator of $\Delta A $ is $\bar{\tau}_{n,n} \tau_{n,n}$ and  $\bar{\tau},\tau$
commute, then, since
$\bar{\tau}^{n+1}_{n,n} \tau^{n+1}_{n,n} =id_{n,n}$, 
we conclude that $(\bar{\tau}_{n,n} \tau_{n,n})^{n+1}=id$. So that $\Delta A$ is a cyclic object. 

Let $k$ be a commutative unital ring. By a paracyclic (resp. cylindrical or cyclic) $k$-module, we mean a paracyclic (resp. cylindrical or cyclic) object in the 
category of $k$-modules.
A parachain complex, by definition, is a graded $k$-module $\mathsf{V}_{\bullet} = ( V_i)_{i \in \mathbb{N}}$\; equipped
with the operators $b: V_i \rightarrow V_{i-1}$ and $ B: V_i \rightarrow V_{i+1} $ such that $ b^2 = B^2 =0 $ and the operator
$T = 1-(bB+Bb)$ is invertible. In the case that $T=1$, the parachain complex is called a mixed complex.
If we denote the graded vector space of
formal power series in the variable $\mathbf{u}$ of degree $-2$ with coefficients in $\mathsf{V}_{\bullet}$ 
by $ \mathsf{V}_{\bullet} \mathbf{\mathit{[[}u\mathit{]]}}$  then we can 
construct the complex $( \mathsf{V}_{\bullet} \mathbf{\mathit{[[}u\mathit{]]}}, b + \mathbf{u} B)$, where $\mathsf{V}_n\mathit{[[}\mathbf{u}\mathit{]]}= 
\prod_{i=0}^{\infty} V_{n+2i}\mathbf{u}^i$. A morphism between the parachain complexes $\mathsf{V}_{\bullet}$ and 
$\mathsf{V'}_{\bullet}$ is a morphism of the complexes $( \mathsf{V}_{\bullet} \mathbf{\mathit{[[}u\mathit{]]}}, b + \mathbf{u} B)$ and
$( \mathsf{V'}_{\bullet} \mathbf{\mathit{[[}u\mathit{]]}}, b' + \mathbf{u} B')$. So a morphism, $ \mathbf{f:} \; \mathsf{V}_{\bullet} \rightarrow
\mathsf{V'}_{\bullet}$ is a sequence of maps $ \mathbf{f}_k : V_i \rightarrow V_{i+2k}, \;  k \ge 0 $ such that
$b'\mathbf{f}_k + B'\mathbf{f}_{k-1} = \mathbf{f}_k b + \mathbf{f}_{k-1} B$. 

Corresponding to any paracyclic module $A$, we can define the parachain complex $\mathsf{C}_{\bullet}(A)$ with the underlying graded 
module   $\mathsf{C}_n(A)= A(n)$ and the operators $ b = \sum_{i=0}^n (-1)^i \partial_i$ and 
$ B = (1 - (-1)^{n+1} \tau ) \sigma N$. Here, $\sigma$ is the extra degeneracy satisfying $\tau \sigma_0 = \sigma \tau $, and
 $ N= \sum_{i=0}^n(-1)^{in} \tau^i $ is the norm operator. So $\mathsf{C}$ is a functor from the category of
paracyclic modules to the category of parachain complexes. Therefore, for any bi-paracyclic module $A$, $Tot(\mathsf{C}(A))$
is a parachain complex with $Tot_n(\mathsf{C}(A))= \sum_{p+q=n} A(p,q)$ and with the operators $ Tot(b) = b + \bar{b} $ and
$Tot(B)=  B + T \bar{B}$, where $ T = 1-(bB+Bb)$.
If we define the normalized chain functor $\mathsf{N}$ from  paracyclic modules to parachain 
complexes with the underlying
graded module $\mathsf{N}_n(A) = A(n)/\sum_{i=0}^n im(s_i) $ and the operators $b ,B$ 
induced from $\mathsf{C}_{\bullet}(A)$, then we have the 
following well-known results (see~\cite{gj93}):\\
1. The quotient map $(\mathsf{C}_{\bullet}(A),b) \rightarrow ( \mathsf{N}_{\bullet}(A) , b)$ is 
a quasi-isomorphism of complexes.\\
2. The Eilenberg-Zilber theorem holds for cylindrical modules, i.e.  \\
For any cylindrical module $A$, there is a natural quasi-isomorphism $\mathbf{f}_0 + \mathbf{u f_1} : 
Tot_{\bullet}(\mathsf{N}(A))
 \rightarrow \mathsf{N}_{\bullet}(\Delta(A))$ of parachain complexes, where $f_0$ is the shuffle map.

When $A$ is a cyclic module then we can construct the related mixed complex $(\mathsf{C}_{\;\bullet}(A)\mathit{[[}\mathbf{u}\mathit{]]} , b + \mathbf{u} B)$.
The homology of the complex $(\mathsf{C}_{\;\bullet}(A) \otimes_{k[\mathbf{u}]} \mathsf{W} , b + \mathbf{u} B)$ 
defines the cyclic homology of $A$ with coefficients in $\mathsf{W}$, where $\mathsf{W}$ is a graded 
$k[\mathbf{u}]$-module with finite homological dimention.
We denote $\mathsf{C}_{\;\bullet}(A) \otimes_{k[\mathbf{u}]} \mathsf{W} $ by the usual notation $ \mathsf{C}(A) \boxtimes
\mathsf{W} $. We know that if
 $\mathsf{W}= k[\mathbf{u}] , k[\mathbf{u , u^{-1}}], k[\mathbf{u , u^{-1}}]/\mathbf{u}k[\mathbf{u}],
k[\mathbf{u}]/\mathbf{u} k[\mathbf{u}]$ we get respectively $ HC_{\bullet}^- (A)$ negative cyclic homology,
$HP_{\bullet}^- (A)$  periodic cyclic homology, $HC_{\bullet}^- (A)$ cyclic homology and $HH_{\bullet}^- (A)$ Hochschild homology~\cite{gj93}.

In this paper the word algebra means an associative, not necessarily commutative, unital algebra over a fixed commutative ring $k$.
Similarly, our Hopf algebras are over $k$ and are not assumed to be commutative or cocommutative. 
The undecorated tensor product $\otimes$ means tensor product over $k$. If $\mathcal{H}$ is a Hopf algebra, we denote its coproduct
by $\Delta : \mathcal{H} \rightarrow \mathcal{H} \otimes \mathcal{H}$, its counit by $\epsilon : \mathcal{H} \rightarrow k$,
its unit by $\eta:k\rightarrow \mathcal{H}$ and its antipode by $S:\mathcal{H} \rightarrow \mathcal{H}$. We will use Sweedler's
notation $\Delta(h)=h^{(0)} \otimes h^{(1)}, \Delta^2(h)=(\Delta \otimes id)\circ \Delta (h)=(id \otimes \Delta )\circ \Delta(h)=h^{(0)} \otimes h^{(1)}\otimes h^{(2)}$, 
etc., where the summation is understood.\\
If $\mathcal{H}$ is a Hopf algebra, the word $\mathcal{H}$-module means a module over the underlying algebra of $\mathcal{H}$. An
algebra $A$ is called a left $\mathcal{H}$-module algebra if $A$ is a left $\mathcal{H}$-module and the multiplication map $A \otimes A \rightarrow A$ and the unit map $k \rightarrow A$ are morphisms of $\mathcal{H}$
-modules. That is $$ h \cdot (ab)=(h^{(0)} \cdot a)(h^{(1)} \cdot b)$$ $$h \cdot 1 = \epsilon(h)1, $$ for all $h \in \mathcal{H}, a,b \in A$. The same definition applies if $\mathcal{H}$
is just a bialgebra. 

Given an $\mathcal{H}$-module algebra $A$, the crossed
product (or smash product) $A \rtimes \mathcal{H}$ of $A$ and $\mathcal{H}$ is, as a $k$-module, $A \otimes \mathcal{H}$, with the
product $$(a \otimes g) (b \otimes h)= a(g^{(0)} \cdot b) \otimes g^{(1)} h.$$
It is an algebra under the above product. We note that the algebra structure of $ A \rtimes \mathcal{H}$ is independent of the
antipode of $\mathcal{H}$ and can be defined when $\mathcal{H}$ is just a bialgebra.
We denote the opposite algebra of $A$ by $A^{op}$. If $\mathcal{H}$ is a bialgebra, $\mathcal{H}^{cop}$ has the same multiplication and
the opposite comultiplication. It is easy to see that if $A$ is a left $\mathcal{H}$-module algebra, then$A^{op}$ is a left $\mathcal{H}^{cop}$-module algebra with the 
same action. 
 It is known that if $\mathcal{H}$ is a Hopf algebra, then $\mathcal{H}^{cop}$ is a Hopf algebra if and 
only if the antipode $S$ is invertible and in that case the antipode of $\mathcal{H}^{cop}$ is $S^{-1}$~\cite{sM95}.\\

%**********************************
\section{The Cylindrical Module $\mathbf {A} \mathbf{\natural} \mathbf{\mathcal{H}}.$}
In this section we introduce the cylindrical module $$A \natural \mathcal{H}= \{  \mathcal{H}^{\otimes (p+1)} \otimes A^{\otimes(q+1)} \}_{p,q \ge 0}, $$ where 
$A$ is an $\mathcal{H}$-module algebra and $\mathcal{H}$ is a Hopf algebra. We define the
operators $\tau_{p,q} ,\partial^{p,q} ,\sigma^{p,q}$ and $\bar{\tau}_{p,q} ,\bar{\partial}^{p,q} ,\bar{\sigma}^{p,q}$ 
as follows:

\begin{eqnarray} \label{eq:opcy1} \notag 
& &\tau_{p,q}( g_0 , \dots , g_p \mid a_{q} , \dots , a_0)=   \\ \notag
& &\hspace{1cm} (g_0^{(0)}, \dots , g_p^{(0)} \mid  a_{q-1}, \dots ,a_{0},S(g_0^{(1)} g_{1}^{(1)} \dots
g_p^{(1)}) \cdot a_q) \\ \notag
& &\partial^{p,q}_i(g_p,\dots,g_0 \mid a_q , \dots , a_0)=  \\ \notag 
& &\hspace{1cm} (g_0,\dots,g_p \mid a_q,\dots , a_{i+1} a_{i},\dots , a_0) \;\;\; 0 \le i <q  \\ 
& &\partial^{p,q}_q (g_0,\dots,g_p \mid a_q , \dots , a_0)=  \\ \notag 
& & \hspace{1cm}(g_0^{(0)},\dots,g_p^{(0)} \mid a_{q-1}, \dots ,a_0(S(g_0^{(1)} g_{1}^{(1)} \dots  g_p^{(1)} ) 
\cdot a_q)) \\ \notag 
& &\sigma^{p,q}_i(g_p,\dots,g_0 \mid a_q , \dots , a_0)= \\ \notag 
& & \hspace{1cm}(g_0,\dots,g_p \mid a_q,\dots , a_{i+1} , 1 , a_{i},\dots , a_0) \;\;\; 0 \le i \le q.  
\end{eqnarray}

\begin{eqnarray} \label{eq:opcy2} \notag
& &\bar{\tau}_{p,q}( g_0 , \dots , g_p \mid a_q , \dots , a_0)=\\ \notag 
& &\hspace{1cm}(g_p^{(0)},g_{0}, \dots , g_{p-1} \mid g_p^{(1)} \cdot a_{q}, \dots
 ,g_p^{(q+1)} \cdot a_0)\\ \notag 
& &\bar{\partial}^{p,q}_i(g_0,\dots,g_p \mid a_q , \dots , a_0)=\\ \notag 
& &\hspace{1cm}(g_0,\dots,g_{i} g_{i+1},\dots,g_p \mid a_q,\dots , a_0) \;\;\; 0 \le i <p  \\ 
& &\bar{\partial}^{p,q}_p (g_0,\dots,g_p \mid a_q , \dots , a_0)=\\ \notag 
& &\hspace{1cm}(g_p^{(0)}g_{0}, \dots , g_{p-1} \mid g_p^{(1)} \cdot a_{q}, \dots
 ,g_p^{(q+1)} \cdot a_0)\\ \notag 
& &\bar{\sigma}^{p,q}_i (g_0,\dots,g_p \mid a_q , \dots , a_0)=\\ \notag 
& &\hspace{1cm}(g_0,\dots,g_{i},1,g_{i+1},\dots,g_p \mid a_q,\dots , a_0) \;\;\; 0 \le i \le p.  
\end{eqnarray}
 
Here we have used Sweedler's notation for the iterated coproducts $\Delta^n(g)=\sum g^{(0)} \otimes g^{(1)} \otimes \dots \otimes g^{(n)}$,
 and we have 
omitted the summation sign for simplicity. 
\begin{theorem}
$A \natural \mathcal{H}$ with the operators defined in (\ref{eq:opcy1}),(\ref{eq:opcy2}) is a cylindrical module.  
\end{theorem}
\begin{proof}
We check only the commutativity of the cyclic operators and the cylindrical condition.

\item[$\bullet \quad$] $\tau_{p,q} \bar{\tau}_{p,q} = \bar{\tau}_{p,q} \tau_{p,q}$
\begin{eqnarray*}
& &(\tau_{p,q} \bar{\tau}_{p,q})(g_0,\dots,g_p \mid a_q,\dots ,a_0)\\
&=& \tau_{p,q}( g_p^{(0)} ,g_{0} \dots , g_{p-1} \mid g_p^{(1)} \cdot a_q ,\dots ,g_p^{(q+1)} \cdot a_0 )\\
&=& ( g_{p}^{(0)} ,g_{0}^{(0)},\dots ,  g_{p-1}^{(0)} \mid  g_{p}^{(3)} \cdot a_{q-1},
\dots ,g_p^{(q+2)} \cdot a_0,(S(g_p^{(1)}g_{0}^{(1)} \dots g_{p-1}^{(1)})\cdot
(g_p^{(2)} \cdot a_q ))\\
&=& ( g_{p}^{(0)} ,g_{0}^{(0)},\dots ,  g_{p-1}^{(0)} \mid  g_{p}^{(1)} \cdot a_{q-1},
\dots ,g_p^{(q)} \cdot a_0,S(g_{0}^{(1)} \dots g_{p-1}^{(1)})
 \cdot a_q )\\
&=& ( g_{p}^{(0)} ,g_{0}^{(0)},\dots ,  g_{p-1}^{(0)} \mid  g_{p}^{(1)} \cdot a_{q-1},
\dots ,g_p^{(q)} \cdot a_0,(g_p^{(q+1)}S(g_{0}^{(1)} \dots g_{p-1}^{(1)}g_p^{(q+2)})) \cdot a_q )\\
&=& \bar{\tau}_{p,q}(g_0^{(0)}, \dots , g_p^{(0)} \mid  a_{q-1},  \dots ,a_0,S(g_0^{(1)}  \dots 
g_p^{(1)}) \cdot a_q)\\
&=& ( \bar{\tau}_{p,q} \tau_{p,q} )(g_0,\dots,g_p \mid a_q,\dots ,a_0).
\end{eqnarray*}

\item[$\bullet \quad$] $\bar{\tau}_{p,q}^{p+1} \;
\tau_{p,q}^{q+1}=\tau_{p,q}^{q+1} \;
\bar{\tau}_{p,q}^{p+1}=id_{p,q}$\\

$\bar{\tau}_{p,q}^{p+1} \tau_{p,q}^{q+1}(g_0,\dots,g_p \mid a_q, \dots , a_0)$
\begin{eqnarray*}
&=& \bar{\tau}_{p,q}^{p+1} \tau_{p,q}^q(g_0^{(0)},\dots , g_p^{(0)} \mid  a_{q-1}, \dots ,a_0, S(g_0^{(1)} g_{1}^{(1)} 
\dots g_p^{(1)}) \cdot a_q)\\ 
& & \qquad \qquad \vdots \\
&=& \bar{\tau}_{p,q}^{p+1}(g_0^{(0)},\dots,g_p^{(0)} \mid S(g_0^{(q+1)} \dots g_p^{(q+1)}) 
\cdot a_q ,S(g_0^{(q)} \dots g_p^{(q)}) \cdot a_{q-1} ,\\
& & \hspace{8.5cm} \dots ,S(g_0^{(1)} \dots g_p^{(1)}) \cdot a_0 )\\
&=& \bar{\tau}_{p,q}^{p}(g_p^{(0)},g_0^{(0)},\dots,g_{p-1}^{(0)} \mid (g_p^{(1)}S(g_0^{(q+1)} \dots g_p^{(2q+2)})) 
\cdot a_q ,
\\& &\hspace{1cm}(g_p^{(2)}S(g_0^{(2q+1)} \dots g_p^{(q)})) \cdot a_{q-1} ,
\dots ,(g_p^{(q+1)}S(g_0^{(1)} \dots g_p^{(q+2)})) \cdot a_0 )\\
&=& \bar{\tau}_{p,q}^{p}(g_0^{(0)},\dots,g_{p-1}^{(0)} \mid S(g_0^{(q+1)} \dots g_{p-1}^{(q+1)}) 
\cdot a_q ,S(g_0^{(q)} \dots g_{p-1}^{(q)}) \cdot a_{q-1} ,\\
& & \hspace{8.5cm} \dots ,S(g_0^{(1)} \dots g_{p-1}^{(1)}) \cdot a_0 )\\
&=&(g_0,\dots,g_p \mid a_q, \dots , a_0).
\end{eqnarray*}
\end{proof}
\begin{corollary}If $S$ is invertible, then $A^{op} \natural \mathcal{H}^{cop}$ is a cylindrical module.
\end{corollary}
To represent the elements of $A^{op} \natural \mathcal{H}^{cop}$ we use direct indexing on $A^{\otimes n}$, 
for example the operators 
$\tau$ and $\bar{\tau}$ will be represented as

\begin{eqnarray*}
& &\tau_{p,q}( g_0 , \dots , g_p \mid a_{0} , \dots , a_q)= \\  
& &\hspace{1cm} (g_0^{(1)}, \dots , g_p^{(1)} \mid  S^{-1}(g_0^{(0)} g_{1}^{(0)} \dots
g_p^{(0)}) \cdot a_q,a_{0}, \dots ,a_{q-1}) \\ 
& &\bar{\tau}_{p,q}( g_0 , \dots , g_p \mid a_0 , \dots , a_q)=\\ 
& &\hspace{1cm}(g_p^{(q+1)},g_{0}, \dots , g_{p-1} \mid g_p^{(0)} \cdot a_0, \dots,
g_p^{(q)} \cdot a_{q})\\  
\end{eqnarray*}
\begin{corollary}For any Hopf algebra $\mathcal{H}$, if $A$ is an $\mathcal{H}$-module algebra, then $\Delta (A \natural \mathcal{H})$ is a cyclic module. If $S$ is invertible then $\Delta (A^{op} \natural \mathcal{H}^{cop})$ is a cyclic module.
\end{corollary}
%******************************
\section{Relation of $\mathbf{\Delta (A \natural \mathcal{H})}$ with the Cyclic Module of the Crossed Product Algebra $\mathbf{ A \rtimes \mathcal{H} }.$} 
We define a map $\phi : \mathsf{C}_{\bullet}(A^{op} \rtimes \mathcal{H}^{cop}) \rightarrow \Delta(A \natural \mathcal{H})$, where 
$\mathsf{C}_{\bullet}(A^{op} \rtimes \mathcal{H}^{cop})$ denotes the cyclic module of the crossed product algebra $A^{op} \rtimes \mathcal{H}^{cop}$, by \\

$\phi ( a_0 \otimes g_0,\dots ,a_n \otimes g_n )=$
\begin{eqnarray*}
& & (g_0^{(0)},g_{1}^{(0)},\dots,g_n^{(0)} \mid S(g_n^{(n+1)}) \cdot a_n,S(g_{n-1}^{(n)}g_{n}^{(n)})\cdot a_{n-1},\dots,\\
& & \hspace{2cm} 
,S(g_{1}^{(2)}g_{2}^{(2)} \dots g_n^{(2)}) \cdot a_1,S(g_0^{(1)}g_{1}^{(1)} \dots g_n^{(1)}) \cdot a_0).   
\end{eqnarray*}
\begin{theorem} $\phi$ defines a cyclic map between $\mathsf{C}_{\bullet}(A^{op} \rtimes \mathcal{H}^{cop})$ and $\Delta (A \natural \mathcal{H})$.
\end{theorem}
\begin{proof}
We show that $\phi$ commutes with the cyclic and simplicial operators:\\
\item[$\bullet \quad$] $ \bar{\tau}_{n,n} \tau_{n,n} \phi = \phi \tau_n^{A^{op} \rtimes \mathcal{H}^{cop}} $ 
\begin{eqnarray*}
& &(\bar{\tau}_{n,n} \tau_{n,n}) \phi (a_0 \otimes g_0, \dots , a_n \otimes g_n)\\
&=&\bar{\tau}_{n,n} \tau_{n,n}(g_0^{(0)},g_{1}^{(0)},\dots,g_n^{(0)} \mid  S(g_n^{(n+2)}) \cdot a_n,S(g_{n-1}^{(n)}
g_{n}^{(n)})\cdot a_{n-1},\dots
\\
& & \hspace{4cm} ,S(g_{1}^{(2)}g_{2}^{(2)} \dots g_n^{(2)}) \cdot a_1,S(g_0^{(1)}g_{1}^{(1)} \dots g_n^{(1)}) 
\cdot a_0)
\\ 
&=& \bar{\tau}_{n,n} (g_0^{(0)},g_{1}^{(0)}, \dots , g_n^{(0)} \mid S(g_{n-1}^{(n+1)}g_{n}^{(n+1)})\cdot a_{n-1},\dots,
S(g_{1}^{(3)}g_{2}^{(3)} \dots g_n^{(3)}) \cdot a_1\\
& &\hspace{3cm},S(g_0^{(2)}g_{1}^{(2)} \dots g_n^{(2)}) 
\cdot a_0,
 (S(g_0^{(1)} g_{1}^{(1)} \dots g_n^{(1)})S(g_n^{n+2}))
\cdot a_n )
\end{eqnarray*}
\begin{eqnarray*}  
&=& (g_{n}^{(0)}, g_0^{(0)}, \dots ,g_{n-1}^{(0)} \mid (g_n^{(1)}S(g_{n-1}^{(n+1)}g_{n}^{(2n+2)}))\cdot a_{n-1},\dots , \\
& &\hspace{1cm}(g_n^{(n+1)}S(g_1^{(3)} \dots g_{n-1}^{(3)} g_{n}^{(n+3)})) \cdot a_1,(g_n^{(n)}S(g_0^{(2)} \dots 
g_{n-1}^{(2)}g_n^{(n+2)})) 
\cdot a_0,\\
& & \hspace{5cm} ,
   (g_n^{(n+1)}S( g_0^{(1)} \dots g_{n-1}^{(1)} g_n^{(n+2)})S(g_n^{2n+3}))
\cdot a_n )\\
&=& (g_n^{(0)}, g_0^{(0)} \dots , g_{n-1}^{(0)}\mid S(g_{n-1}^{(n+1)})\cdot a_{n-1} ,\dots,\\
& &\hspace{2cm}S( g_1^{(3)}\dots g_{n-1}^{(3)}) \cdot a_1,S(g_0^{(2)} \dots g_{n-1}^{(2)}) 
\cdot a_0,S( g_n^{(1)}g_0^{(1)} \dots g_{n-1}^{(1)})
\cdot a_n )
\\ 
&=& \phi (a_n \otimes g_n,a_0 \otimes g_0 , \dots , a_{n-1} \otimes g_{n-1})
= \phi (\tau_n^{A^{op} \rtimes \mathcal{H}^{cop}}( a_0 \otimes g_0, \dots , a_n \otimes g_n)).\\
\end{eqnarray*}

\item[$\bullet \quad$] $ \bar{\partial}_n^{n,n} \partial_n^{n,n} \phi = \phi \partial_n^{A^{op} \rtimes \mathcal{H}^{cop}} $ 
\begin{eqnarray*}
& & \bar{\partial}_n^{n,n} \partial_n^{n,n} \phi (a_0 \otimes g_0,\dots , a_n \otimes g_n)\\
&=&\bar{\partial}_n^{n,n} \partial_n^{n,n}(g_0^{(0)},g_{1}^{(0)},\dots,g_n^{(0)} \mid  S(g_n^{(n+2)}) \cdot a_n,S(g_{n-1}^{(n)}
g_{n}^{(n)})\cdot a_{n-1},\dots
\\
& & \hspace{4cm} ,S(g_{1}^{(2)}g_{2}^{(2)} \dots g_n^{(2)}) \cdot a_1,S(g_0^{(1)}g_{1}^{(1)} \dots g_n^{(1)}) 
\cdot a_0)
\\ 
&=& \bar{\partial}_n^{n,n} (g_0^{(0)},g_{1}^{(0)}, \dots , g_n^{(0)} \mid S(g_{n-1}^{(n+1)}g_{n}^{(n+1)})\cdot a_{n-1},\dots,
S(g_{1}^{(3)}g_{2}^{(3)} \dots g_n^{(3)}) \cdot a_1\\
& &\hspace{3cm},(S(g_0^{(2)}g_{1}^{(2)} \dots g_n^{(2)}) 
\cdot a_0)
 ((S(g_0^{(1)} g_{1}^{(1)} \dots g_n^{(1)})S(g_n^{n+2}))
\cdot a_n ))\\ 
&=& (g_{n}^{(0)}g_0^{(0)}, \dots ,g_{n-1}^{(0)} \mid (g_n^{(1)}S(g_{n-1}^{(n+1)}g_{n}^{(2n+2)}))\cdot a_{n-1},\dots , \\
& &\hspace{6cm}(g_n^{(n+1)}S(g_1^{(3)} \dots g_{n-1}^{(3)} g_{n}^{(n+3)})) \cdot a_1,\\
& &  ,((g_n^{(n)}S(g_0^{(2)} \dots 
g_{n-1}^{(2)}g_n^{(n+2)})) 
\cdot a_0)
   ((g_n^{(n+1)}S( g_0^{(1)} \dots g_{n-1}^{(1)} g_n^{(n+2)})S(g_n^{2n+3}))
\cdot a_n ))\\ 
&=& (g_n^{(0)} g_0^{(0)} \dots , g_{n-1}^{(0)}\mid S(g_{n-1}^{(n+1)})\cdot a_{n-1} ,\dots,\\
& &\hspace{2cm}S( g_1^{(3)}\dots g_{n-1}^{(3)}) \cdot a_1,(S(g_0^{(2)} \dots g_{n-1}^{(2)}) 
\cdot a_0)(S( g_n^{(1)}g_0^{(1)} \dots g_{n-1}^{(1)})
\cdot a_n) )
\\
&=& (g_n^{(0)} g_0^{(0)} \dots , g_{n-1}^{(0)}\mid S(g_{n-1}^{(n+1)})\cdot a_{n-1} ,\dots,\\
& &\hspace{2cm}S( g_1^{(3)}\dots g_{n-1}^{(3)}) \cdot a_1,S( g_n^{(1)}g_0^{(1)} \dots g_{n-1}^{(1)})
\cdot((g_n^{(2)} \cdot a_0) a_n) )
\\  
&=& \phi ((a_n \otimes g_n)(a_0 \otimes g_0) , \dots , a_{n-1} \otimes g_{n-1})\\
&=& \phi \partial_n^{A^{op} \rtimes \mathcal{H}^{cop}} (a_0 \otimes g_0 , \dots , a_n \otimes g_n).
\end{eqnarray*}
With the same argument one can check that
$ \bar{\partial}_i^{n,n} \partial_i^{n,n} \phi = \phi \partial_i^{A^{op} \rtimes \mathcal{H}^{cop}} $
where $ 0 \le i < n $ and $\bar{\sigma}_i^{n,n} \sigma_i^{n,n} \phi = \phi \sigma_i^{A^{op} \rtimes \mathcal{H}^{cop}}. $\\ 
\end{proof}
\begin{theorem} \label{th:diag} 
Let $\mathcal{H}$ be a Hopf algebra and $A$ an $\mathcal{H}$-module algebra. Then we have an isomorphism of cyclic modules
$$\Delta (A \natural \mathcal{H}) \cong \mathsf{C}_{\bullet}(A^{op} \rtimes \mathcal{H}^{cop}).$$
If $S$ is invertible then $$\Delta (A^{op} \natural \mathcal{H}^{cop}) \cong \mathsf{C}_{\bullet}(A \rtimes \mathcal{H}).$$
\end{theorem}
\begin{proof}
One can check that  $\psi$ is a cyclic map and $\phi \circ \psi = \psi \circ \phi =id$ , where $\psi$ is defined by\\

$\psi (g_0, \dots , g_n \mid a_n , \dots ,a_0)=$  
\begin{eqnarray*}
(g_0^{(1)} \dots g_{n-1}^{(n)} g_n^{(n+1)} ) \cdot a_0 \otimes g_0^{(0)} ,( g_1^{(1)} \dots g_{n-1}^{(n-1)}g_n^{(n)}) \cdot a_1 \otimes g_1^{(0)},
\dots , (g_n^{(1)}) \cdot a_n \otimes g_n^{(0)}).
\end{eqnarray*}
Precisely,\\

$\phi \circ \psi (g_n, \dots , g_0 \mid a_n , \dots ,a_0)=$
\begin{eqnarray*}
& &(g_0^{(0)},g_{1}^{(0)},\dots,g_n^{(0)} \mid (S(g_{n}^{(n+1)}) (g_n^{(n+2)})) \cdot a_{n},
(S(g_{n-1}^{(n)} g_{n}^{(n)})(g_{n-1}^{(n+1)}g_{n}^{(n+3)}))\cdot a_{n-1}, \\
& &\hspace{2cm}\dots,(S(g_1^{(2)} g_{2}^{(2)} \dots g_n^{(2)})(g_1^{(3)} g_{2}^{(5)}\dots  g_n^{(2n+1)})) \cdot a_1,\\
& &\hspace{4cm}(S(g_0^{(1)}  \dots g_{n-1}^{(1)} g_n^{(1)} )(g_0^{(2)}g_1^{(4)} \dots g_n^{(2n+2)})) \cdot a_0)\\
&=&(g_0, \dots , g_n \mid a_n , \dots ,a_0),
\end{eqnarray*}
and\\

$\psi \circ \phi ( a_0 \otimes g_0,\dots ,a_n \otimes g_n )=$
\begin{eqnarray*}
& &(( g_0^{(1)}g_{1}^{(2)} \dots g_n^{(n+1)})S( g_0^{(2)}g_{1}^{(3)} \dots g_n^{(n+2)} )) \cdot a_0 \otimes g_0^{(0)},\\
& & \hspace{2cm}(( g_1^{(1)} g_{2}^{(2)} \dots g_n^{(n)})S( g_1^{(4)}  \dots g_{n-1}^{(n+2)} g_n^{(n+3)})) \cdot a_1 \otimes g_1^{(0)},\dots,\\
& &\hspace{.5cm}((g_{n-1}^{(1)}g_{n}^{(2)} )S(g_{n-1}^{(2n)}g_{n}^{(2n+1)} )) \cdot a_{n-1} \otimes g_{n-1}^{(0)},
(g_{n}^{(1)}S(g_n^{(2n+2)})) \cdot a_{n} \otimes g_{n}^{(0)})\\
&=& (a_0 \otimes g_0,\dots ,a_n \otimes g_n).   
\end{eqnarray*}

If $S$ is invertible then we can see that $\phi$ and $\psi$ will be represented as\\

$\phi ( a_0 \otimes g_0,\dots ,a_n \otimes g_n )=$
\begin{eqnarray*}
& & (g_0^{(1)},g_1^{(2)},\dots,g_n^{(n+1)} \mid  S^{-1}(g_0^{(0)}g_1^{(1)} \dots g_n^{(n)}) \cdot a_0,S^{-1}(g_1^{(0)}g_2^{(1)} \dots g_n^{(n-1)}) \cdot a_1,\dots \\
& & \hspace{5cm} ,
S^{-1}(g_{n-1}^{(0)}g_n^{(1)})\cdot a_{n-1},S^{-1}(g_n^{(0)}) \cdot a_n),   
\end{eqnarray*}

$\psi (g_0, \dots , g_n \mid a_0 , \dots ,a_n)=$  
\begin{eqnarray*}
((g_0^{(0)} g_1^{(0)} \dots g_n^{(0)}) \cdot a_0 \otimes g_0^{(1)} ,(g_1^{(1)} \dots g_n^{(1)}) \cdot a_1 \otimes g_1^{(2)},
\dots , g_n^{(n)} \cdot a_n \otimes g_n^{(n+1)}).
\end{eqnarray*}

\end{proof}
%********************************
\section{Cyclic Homology of $\mathsf{C}_{\bullet} (A \rtimes \mathcal{H}).$}
In this section we give a spectral sequence to compute the cyclic homology of $\mathsf{C}_{\bullet}(A^{op} \rtimes \mathcal{H}^{cop})$
and $\mathsf{C}_{\bullet}(A \rtimes \mathcal{H})$ when $S$ is invertible.
 
By using the Eilenberg-Zilber Theorem for cylindrical modules, combined 
with Theorem~\ref{th:diag} we get the following result.
\begin{theorem} \label{th:Apei-zi}
If $A$ is an $\mathcal{H}$-module algebra, then there is an isomorphism of cyclic homology groups 
\[ HC_{\bullet}(
\text{Tot} (A \natural \mathcal{H});\mathsf{W}) \cong HC_{\bullet} ( \Delta (A \natural \mathcal{H} );\mathsf{W}) \cong
HC_{\bullet}(A^{op} \rtimes \mathcal{H}^{cop};\mathsf{W}). \] 
If $S$ is invertible, then we have 
\[ HC_{\bullet}(
\text{Tot} (A^{op} \natural \mathcal{H}^{cop});\mathsf{W}) \cong HC_{\bullet} ( \Delta (A^{op} \natural \mathcal{H}^{cop} );\mathsf{W}) \cong
HC_{\bullet}(A \rtimes \mathcal{H};\mathsf{W}). \]   
\end{theorem}
Now we define an action of $\mathcal{H}$ on the first column of $A \natural \mathcal{H}$, denoted $(A^{op})^{\natural}_{\mathcal{H}^{cop}}$ 
= $ \{ \mathcal{H} \otimes A^{\otimes (n+1)} \}_{n \ge 0}$, by 
\begin{equation} \label{eq:ac}
h \cdot (g \mid a_n, \dots , a_0 ) = (h^{(1)}  g S(h^{(0)}) \mid h^{(2)} \cdot a_n , \dots , h^{(n+2)} \cdot a_0). 
\end{equation}
If $S$ is invertible, we denote the first column of $A^{op} \natural \mathcal{H}^{cop}$ by $A^{\natural}_{\mathcal{H}}$= 
$ \{ \mathcal{H} \otimes A^{\otimes (n+1)} \}_{n \ge 0}$ and we define an action on it by
\begin{equation}\label{eq:ac1} 
h \cdot (g \mid a_0, \dots , a_n ) = (h^{(n+2)}  g S(h^{(0)}) \mid h^{(1)} \cdot a_0 , \dots , h^{(n+1)} \cdot a_n). 
\end{equation}

We define $C_{\bullet}^{\mathcal{H}^{cop}}(A^{op})$ as the co-invariant space of $(A^{op})^{\natural}_{\mathcal{H}^{cop}}$  under the 
action defined in~(\ref{eq:ac}), i.e.,
\[ C^{\mathcal{H}^{cop}}_{\bullet} (A^{op})=(A^{op})^{\natural}_{\mathcal{H}^{cop}}/ span \{ h \cdot x - \epsilon(h) x \;|\; 
h \in \mathcal{H}^{cop}, x \in (A^{op})^{\natural}_{\mathcal{H}^{cop}} \}. \]

Similarly we define $C_{\bullet}^{\mathcal{H}}(A)$ as the co-invariant space of $A^{\natural}_{\mathcal{H}}$ under 
the action defined in~(\ref{eq:ac1}). It is easy to see that when $S$ is invertible the equivalent co-invariant space 
defined under the action~(\ref{eq:ac}) is
the same as the one is defined by~(\ref{eq:ac1}). \\
We define the following operators on $C_{\bullet}^{\mathcal{H}^{cop}}(A^{op})$, 
\begin {eqnarray}\label{eq:equ5} \notag
& &\tau_n(g \mid a_n , \dots , a_0 )= ( g^{(0)} \mid a_{n-1},\dots,a_0,S(g^{(1)}) \cdot a_n),  \\ \notag
& &\partial_i(g \mid a_n , \dots , a_0 )= ( g \mid a_n,\dots ,a_{i+1} a_{i},\dots,a_{0}),  \;\;\;0 \le i<n, \\
& &\partial_n(g \mid a_n , \dots , a_0 )= ( g^{(0)} \mid a_{n-1},\dots,a_0(S(g^{(1)}) \cdot a_n)),\\ \notag
& &\sigma_i(g \mid a_n , \dots , a_0 )= ( g \mid a_n,\dots ,a_{i+1},1,a_{i},\dots,a_{0}), \;\;\; 0 \le i \le n. \notag
\end{eqnarray}

One can check that these operators are well defined on the co-invariant space. For example we see that
\begin{eqnarray*}
& & \tau_n(h \cdot (g \mid a_n, \dots , a_0 ) - \epsilon (h) (g \mid a_n, \dots , a_0 )) \\
&=& h^{(1)} \cdot (g^{(0)} \mid a_{n-1}, \dots, a_0,   S(g^{(1)}S(h^{(0)})h^{(2)}) \cdot a_n)- \\
& & \hspace{4cm} \epsilon (h^{(1)})( g^{(0)} \mid a_{n-1},  \dots , a_0,S(g^{(1)}S(h^{(0)})h^{(2)}) \cdot a_n).  
\end{eqnarray*}
 
\begin{theorem}
$C_{\bullet}^{\mathcal{H}^{cop}}(A^{op})$ with the operators defined in $(\ref{eq:equ5})$ is a cyclic module.
\end{theorem}
\begin{proof}
We only check that $\tau_n^{n+1}=id$, the other identities are similar to check. We see that
\begin{eqnarray*}
& &(g^{(0)} \mid g^{(1)} \cdot a_n, \dots , g^{(n+1)} \cdot a_0) = 
(g^{(2)} g^{(0)} S(g^{(1)}) \mid g^{(3)} \cdot a_n, \dots , g^{(n+3)} \cdot a_0)\\
&=&g^{(1)} \cdot (g^{(0)} \mid a_n, \dots , a_0) \equiv \epsilon(g^{(1)}) (g^{(0)} \mid a_n, \dots , a_0) = 
(g \mid  a_n, \dots , a_0).
\end{eqnarray*}
Therefore we have,
\begin{eqnarray*}
& &\tau_n^{n+1}(g \mid a_n , \dots , a_0 )= \tau_n^{n+1}(g^{(0)} \mid g^{(1)} \cdot a_n, \dots , g^{(n+1)} \cdot a_0)\\
&=&\tau_n^{n}(g^{(0)} \mid g^{(3)} \cdot a_{n-1}, \dots , g^{(n+2)} \cdot a_0, (S(g^{(1)})g^{(2)}) \cdot a_n)\\
&=&\tau_n^{n}(g^{(0)} \mid g^{(1)} \cdot a_{n-1}, \dots , g^{(n)} \cdot a_0, a_n)\\
& &\vdots\\
&=&(g \mid a_n , \dots , a_0 ).
\end{eqnarray*}
\end{proof}
\begin{lem}
If $S$ is invertible then  
$C_{\bullet}^{\mathcal{H}}(A)$ with the operators defined as follows is a cyclic module.
\begin {eqnarray*}
& &\tau_n(g \mid a_0 , \dots , a_n )= ( g^{(1)} \mid S^{-1}(g^{(0)}) \cdot a_n,a_0,\dots,a_{n-1}),  \\ \notag
& &\partial_i(g \mid a_0 , \dots , a_n )= ( g \mid a_0,\dots ,a_{i} a_{i+1},\dots,a_{n}),  \;\;\;0 \le i<n, \\
& &\partial_n(g \mid a_0 , \dots , a_n )= ( g^{(1)} \mid (S^{-1}(g^{(0)}) \cdot a_n)a_{0},\dots,a_{n-1}),\\ \notag
& &\sigma_i(g \mid a_0 , \dots , a_n )= ( g \mid a_0,\dots ,a_{i},1,a_{i+1},\dots,a_{n}), \;\;\; 0 \le i \le n. \notag
\end{eqnarray*}
\end{lem}
\begin{proof}
We can see that these operators are directly well defined on the equivariant space, i.e., they
commute with the action defined in~(\ref{eq:ac1}). 
\end{proof}
Let $M$ be a left $\mathcal{H}$-module. We denote 
the space of $p$-chains on $\mathcal{H}$ with values in $M$ by $C_p(\mathcal{H},M)=\mathcal{H}^{\otimes p} \otimes M$, and we define the boundary $\delta : C_p(\mathcal{H},M) \rightarrow
C_{p-1}(\mathcal{H},M)$ by 
$$\delta (g_1,g_2,\dots , g_p,m)= $$
\begin{eqnarray*} 
& & \epsilon (g_1) (g_2,\dots ,g_p,m) +\\
& &\sum_{i=1}^{p-1} (-1)^i (g_1 , \dots ,g_i g_{i+1},\dots , g_p ,m) + (-1)^p (g_1,\dots,g_{p-1},g_p \cdot m). 
\end{eqnarray*}
We denote the $p$-th homology of the complex $(C_{\bullet}(\mathcal{H},M),\delta)$ by $H_p(\mathcal{H},M)$.\\
Let $\mathsf{C}_q((A^{op})_{\mathcal{H}^{cop}}^{\natural})= 
\mathcal{H} \otimes A^{\otimes(q+1)}.$ Since it is an $\mathcal{H}$-module via~(\ref{eq:ac}), we can construct
$H_p(\mathcal{H},\mathsf{C}_q((A^{op})^{\natural}_{\mathcal{H}^{cop}}))$.
        
Now we replace the complex $(A \natural \mathcal{H},(\partial,\sigma,\tau),(\bar{\partial},\bar{\sigma},\bar{\tau}))$ 
with an isomorphic complex 
$$(\mathsf{C}(\mathcal{H},\mathsf{C} ( (A^{op})^{\natural}_{\mathcal{H}^{cop}})),( \mathfrak{d},\mathfrak{s},\mathfrak{t}),(\bar{\mathfrak{d}},\bar{\mathfrak{s}},\bar{\mathfrak{t}})),$$
under the transformation defined by 
the maps $\beta : (A \natural \mathcal{H})_{p,q} \rightarrow \mathsf{C}_p(\mathcal{H},\mathsf{C}_q((A^{op})_{\mathcal{H}^{cop}}^{\natural}))$
and $\gamma :  \mathsf{C}_p(\mathcal{H},\mathsf{C}_q((A^{op})_{\mathcal{H}^{cop}}^{\natural})) \rightarrow  (A \natural \mathcal{H})_{p,q}$
\begin{eqnarray*}
& &\beta ( g_0, \dots , g_p \mid a_q , \dots , a_0 ) = ( g_1^{(1)}, \dots , g_p^{(1)} \mid g_0 g_1^{(0)} \dots g_p^{(0)} \mid
a_q, \dots , a_0 ) \\
& &\gamma ( g_1 , \dots , g_p \mid g \mid a_q , \dots , a_0) = (g S(g_1^{(0)} \dots g_p^{(0)}),g_1^{(1)}, \dots , g_p^{(1)}  \mid 
a_q, \dots , a_0 ).
\end{eqnarray*}

We see that $\beta \circ \gamma = \gamma \circ \beta = id $. We find the operators,
$( \mathfrak{d},\mathfrak{s},\mathfrak{t}),(\bar{\mathfrak{d}},\bar{\mathfrak{s}},\bar{\mathfrak{t}})$
under this transformation. \\

$ \hspace{0.2cm} \bullet \quad $ 
First we compute $\bar{\mathfrak{b}} = \beta \bar{b} \gamma$. Since 
\begin{eqnarray*}
\bar{b} ( g_0, \dots , g_p \mid a_q , \dots , a_0 ) &=& \sum_{0 \le i < p} (-1)^i (g_0 , \dots , g_i g_{i+1} ,\dots ,
 g_p \mid a_q , \dots , a_0 ) \\
&+& (-1)^p (g_p^{(0)}g_0 , g_1, \dots , g_{p-1} \mid g_p^{(1)} \cdot a_q , \dots , g_p^{(q+1)} \cdot a_0), 
\end{eqnarray*}
we see that,
\begin{eqnarray*}
& &\bar{\mathfrak{b}}( g_1 , \dots , g_p \mid g \mid a_q , \dots , a_0) = \beta \bar{b} (g S(g_1^{(0)} \dots g_p^{(0)}),
g_1^{(1)}, \dots g_p^{(1)} \mid a_q , \dots , a_0 )\\
&=& \epsilon(g_1)( g_2, \dots , g_p \mid  g  \mid a_q , \dots , a_0)\\ 
&+& \sum _{0 \le i < p} (-1)^i ( g_1,\dots , g_i g_{i+1} , \dots , g_p  \mid g  \mid a_q , \dots , a_0)\\ 
&+& (-1)^p ( g_1, \dots , g_{p-1} \mid g_p^{(1)} g S(g_p^{(0)}) \mid g_p^{(2)} \cdot a_q , \dots , g
_p^{(q+2)} \cdot a_0)\\
&=& \delta  (g_0, \dots , g_p \mid a_q , \dots , a_0 ).
\end{eqnarray*}
So we have,
\begin{eqnarray*}
& &\bar{\mathfrak{d}}_0 (g_1, \dots ,g_p \mid g \mid a_q, \dots , a_0) = \epsilon(g_1)(g_2 , \dots , g_p \mid g \mid a_q , \dots , a_0)\\
& &\bar{\mathfrak{d}}_i (g_1, \dots ,g_p \mid g \mid a_q, \dots , a_0) = (g_1 , \dots , g_i g_{i+1} ,\dots , g_p \mid g \mid a_q , \dots , a_0)\\
& &\bar{\mathfrak{d}}_p (g_1, \dots ,g_p \mid g \mid a_q, \dots , a_0) = (g_1 , \dots , g_{p-1} \mid g_p \cdot 
( g \mid a_q , \dots ,  a_0)).
\end{eqnarray*}
When $S$ is invertible we can see that
\begin{eqnarray*}
& &\bar{\mathfrak{b}}( g_1 , \dots , g_p \mid g \mid a_q , \dots , a_0) = \beta \bar{b} (g S(g_1^{(0)} \dots g_p^{(0)}),
g_1^{(1)}, \dots g_p^{(1)} \mid a_0 , \dots , a_q )\\
&=& \epsilon(g_1)( g_2, \dots , g_p \mid  g  \mid a_0 , \dots , a_q)\\ 
&+& \sum _{0 \le i < p} (-1)^i ( g_1,\dots , g_i g_{i+1} , \dots , g_p  \mid g  \mid a_0 , \dots , a_q)\\ 
&+& (-1)^p ( g_1, \dots , g_{p-1} \mid g_p^{(q+2)} g S(g_p^{(0)}) \mid g_p^{(1)} \cdot a_0 , \dots , g
_p^{(q+1)} \cdot a_q)\\
&=& \delta  (g_0, \dots , g_p \mid a_q , \dots , a_0 ).
\end{eqnarray*}
$ \hspace{0.2cm} \bullet \quad $ Next we compute  $\bar{\mathfrak{s}}_i = \beta \bar{\sigma}_i \gamma$ and
$\bar{\mathfrak{t}} = \beta \bar{\tau} \gamma$. The result is 
\begin{eqnarray*}
& &\bar{\mathfrak{s}}_i (g_1,\dots , g_p \mid a_q , \dots , a_0) 
= ( g_1, \dots ,g_i,1 , g_{i+1}, \dots  g_p \mid g  \mid a_q , \dots , a_0)\\
& & \bar{\mathfrak{t}}(g_1,\dots,g_p \mid g \mid a_q, \dots ,a_0)=\\
& & ( g^{(1)}S(g_1^{(0)} \dots g_{p-1}^{(0)}g_p^{(0)}),g_1^{(2)}, \dots ,g_{p-1}^{(2)} \mid 
g_p^{(2)}  g^{(0)} S(g_p^{(1)})  \mid g_p^{(3)} \cdot a_q, \dots ,  g_p^{(q+3)} \cdot a_0)\\
& & \hspace{1cm}=( g^{(1)}S(g_1^{(0)} \dots g_{p-1}^{(0)}g_p^{(0)}),g_1^{(2)}, \dots ,g_{p-1}^{(2)} \mid 
g_p^{(1)} \cdot (g^{(0)} \mid a_q, \dots ,  a_0)).
\end{eqnarray*}
When $S$ is invertible we have
\begin{eqnarray*}
& &\bar{\mathfrak{s}}_i (g_1,\dots , g_p \mid a_0 , \dots , a_q) 
= ( g_1, \dots ,g_i,1 , g_{i+1}, \dots  g_p \mid g  \mid a_0 , \dots , a_q)\\
& & \bar{\mathfrak{t}}(g_1,\dots,g_p \mid g \mid a_q, \dots ,a_0)= \\
& & ( g^{(1)}S(g_1^{(0)} \dots g_{p-1}^{(0)}g_p^{(0)}),g_1^{(2)}, \dots ,g_{p-1}^{(2)} \mid 
g_p^{(q+3)}  g^{(0)} S(g_p^{(1)})  \mid g_p^{(2)} \cdot a_0, \dots ,  g_p^{(q+2)} \cdot a_q)\\
& &\hspace{1cm} ( g^{(1)}S(g_1^{(0)} \dots g_{p-1}^{(0)}g_p^{(0)}),g_1^{(2)}, \dots ,g_{p-1}^{(2)} \mid 
g_p^{(1)} \cdot (g^{(0)} \mid a_0, \dots ,  a_q)).
\end{eqnarray*}
$ \hspace{0.2cm} \bullet \quad $ Now we compute the operator  $\mathfrak{b} = \beta b \gamma$. We have 
\begin{eqnarray*}
& &\mathfrak{b}(g_1, \dots , g_p \mid g \mid a_q , \dots , a_0) = \beta b (g S(g_1^{(0)} \dots g_p^{(0)}),
g_1^{(1)}, \dots ,g_p^{(1)} \mid a_q , \dots , a_0 )\\
&=&  \sum_{ 0 \le i <q} (-1)^i( 
g_1, \dots, g_p \mid g  \mid a_q , \dots ,a_{i+1} a_{i}, \dots , a_0 )\\  
&+& (-1)^q (g_1^{(1)}, \dots , g_p^{(1)} \mid g^{(0)} \mid a_{q-1} , \dots , a_0( S(g^{(1)} S(g_1^{(0)} \dots g_p^{(0)}) g_1^{(2)} \dots g_p^{(2)} ) \cdot a_q)). 
\end{eqnarray*}
So we conclude that,
\begin{eqnarray*}
& &\mathfrak{d}_i(g_1,\dots,g_p \mid g \mid a_q,\dots,a_0)=(g_1 , \dots , g_p \mid g \mid a_q, \dots,a_{i+1} a_{i},\dots,a_{0}).\\
& &\mathfrak{d}_q(g_1,\dots,g_p \mid g \mid a_q,\dots,a_0)=\\
& &\hspace{1cm} (g_1^{(1)}, \dots , g_p^{(1)} \mid g^{(0)} \mid a_{q-1} , \dots , a_0( S(g^{(1)} S(g_1^{(0)} \dots g_p^{(0)}) g_1^{(2)} \dots g_p^{(2)} ) \cdot a_q)). 
\end{eqnarray*}
When  $S$ is invertible we conclude that
\begin{eqnarray*}
& &\mathfrak{b}(g_1, \dots , g_p \mid g \mid a_q , \dots , a_0) = \beta b (g S(g_1^{(0)} \dots g_p^{(0)}),
g_1^{(1)}, \dots ,g_p^{(1)} \mid a_0 , \dots , a_q )\\
&=&  \sum_{ 0 \le i <q} (-1)^i( 
g_1, \dots, g_p \mid g  \mid a_0 , \dots ,a_{i} a_{i+1}, \dots , a_q )\\  
&+& (-1)^q (g_1, \dots , g_p \mid g^{(1)} \mid  (S^{-1}(g^{(0)}) \cdot a_q)a_0, \dots, a_{q-1} ). 
\end{eqnarray*}

$ \hspace{0.2cm} \bullet \quad $ We consider  $\mathfrak{s}_i = \beta \sigma_i \gamma$ and $\mathfrak{t} = 
\beta \tau \gamma$. We obtain 
\begin{eqnarray*}
& &\mathfrak{s}_i (g_1,\dots,g_p \mid g \mid a_q, \dots ,a_0)= 
(g_1, \dots , g_p \mid g  \mid  a_q , \dots ,a_{i+1},1,a_{i},\dots, a_0).
\end{eqnarray*}

Also,
\begin{eqnarray*}
& & \mathfrak{t} (g_1, \dots , g_p \mid g \mid a_q. \dots , a_0 ) = \\
& & (g_1^{(1)} , \dots ,g_p^{(1)} \mid g^{(0)} \mid a_{q-1},\dots,a_0, S(g^{(1)} S(g_1^{(0)} \dots g_p^{(0)} ) 
g_1^{(2)} \dots g_p^{(2)}) \cdot a_q ). 
\end{eqnarray*}

When $S$ is invertible these computations will be as
\begin{eqnarray*}
& &\mathfrak{s}_i (g_1,\dots,g_p \mid g \mid a_q, \dots ,a_0)= 
(g_1, \dots , g_p \mid g  \mid  a_0 , \dots 
,a_{i},1,a_{i+1},\dots, a_q),\\
& & \mathfrak{t} (g_1, \dots , g_p \mid g \mid a_q. \dots , a_0 ) = 
(g_1 , \dots ,g_p \mid g^{(1)} \mid  S^{-1}(g^{(0)}) \cdot a_q ,a_0,\dots,  a_{q-1}). 
\end{eqnarray*}

By the above computations we can state the following theorems.
\begin{theorem} \label{th:iso1}
The complex $(\mathsf{C} (A \natural \mathcal{H}) \boxtimes \mathsf{W}, b + \mathbf{u} B , \bar{b} + \mathbf{u} \bar{B}) $ 
is isomorphic to the complex 
$(\mathsf{C}(\mathcal{H}, \mathsf{C}((A^{op})^{\natural}_{\mathcal{H}^{cop}} \boxtimes \mathsf{W})), \mathfrak{b} + \mathbf{u} \mathfrak{B} , \bar{\mathfrak{b}} + \mathbf{u} \bar{\mathfrak{B}}) $
where $\mathfrak{\bar{b}}$ is the Hopf-module boundary. 
\end{theorem}
\begin{theorem} \label{th:iso2}
For  $p \ge 0$, $H_p(\mathcal{H},\mathsf{C}_{\bullet} ((A^{op})^{\natural}_{\mathcal{H}^{cop}}))$, 
with the operators defined as follows,
are cyclic modules, where
$H_0(\mathcal{H},\mathsf{C}_{\bullet}((A^{op})^{\natural}_{\mathcal{H}^{cop}})) 
\cong C_{\bullet}^{\mathcal{H}^{cop}}(A^{op}) $, 
\begin {eqnarray*} \label{eq:equ}
& &\mathfrak{t} (g_1, \dots ,g_p \mid g \mid a_q , \dots , a_0 )=\\
& &(g_1^{(1)} , \dots ,g_p^{(1)} \mid g^{(0)} \mid a_{q-1},  \dots , a_0, S(g^{(1)} S(g_1^{(0)} \dots g_p^{(0)} ) 
g_1^{(2)} \dots g_p^{(2)}) \cdot a_q),\\ 
& &\mathfrak{d}_i(g_1,\dots,g_p \mid g \mid a_q,\dots,a_0)=(g_1 , \dots , g_p \mid g \mid a_q, 
\dots,a_{i+1} a_{i},\dots,a_{0}),\\
& &\mathfrak{d}_q(g_1,\dots,g_p \mid g \mid a_q,\dots,a_0)=\\
& &\hspace{2cm} (g_1^{(1)} , \dots ,g_p^{(1)} \mid g^{(0)} \mid a_{q-1},  \dots , a_0( S(g^{(1)} S(g_1^{(0)} 
\dots g_p^{(0)} ) g_1^{(2)} \dots g_p^{(2)}) \cdot a_q)),\\ 
& &\mathfrak{s}_i(g_1, \dots ,g_p \mid g \mid a_q , \dots , a_0 )=(g_1, \dots , g_p \mid g \mid  a_q , 
\dots ,a_{i+1},1,a_{i},\dots, a_0).
\end{eqnarray*} 
\end{theorem}

\begin{prop} 
If $S$ is invertible then for  $p \ge 0$, $H_p(\mathcal{H},\mathsf{C}_{\bullet} (A^{\natural}_{\mathcal{H}}))$, with the operators defined as follows,
are cyclic modules, where
$H_0(\mathcal{H},\mathsf{C}_{\bullet}(A^{\natural}_{\mathcal{H}})) \cong C_{\bullet}^{\mathcal{H}}(A) $: 
\begin {eqnarray*} \label{eq:equ}
& &\mathfrak{t} (g_1, \dots ,g_p \mid g \mid a_0 , \dots , a_q )=
(g_1 , \dots ,g_p \mid g^{(1)} \mid S^{-1}(g^{(0)}) \cdot a_q,a_{0},  \dots , a_{q-1} ),\\ 
& &\mathfrak{d}_i(g_1,\dots,g_p \mid g \mid a_0,\dots,a_q)=(g_1 , \dots , g_p \mid g \mid a_0, \dots,a_{i} a_{i+1},
\dots,a_{q}),\\
& &\mathfrak{d}_q(g_1,\dots,g_p \mid g \mid a_0,\dots,a_q)=
(g_1 , \dots ,g_p \mid g^{(1)} \mid (S^{-1}(g^{(0)}) \cdot a_q)a_{0},  \dots , a_{q-1} ),\\
& &\mathfrak{s}_i(g_1, \dots ,g_p \mid g \mid a_0 , \dots , a_q )=(g_1, \dots , g_p \mid g \mid  a_0 , 
\dots ,a_{i},1,a_{i+1},\dots, a_q).
\end{eqnarray*}
\end{prop}

To compute the homology of the mixed complex $(Tot(A \natural \mathcal{H}), b + \bar{b} + \mathbf{u} (B + \bar{B})),$
we filter it by the subspaces \[
\mathsf{F}_i Tot_n((A \natural \mathcal{H} ) \boxtimes \mathsf{W}) = \sum_{q \le i,\;p+q=n} (\mathcal{H}^{\otimes(p+1)} \otimes A^{\otimes(q+1)}) \boxtimes \mathsf{W}. \]
If we separate the operator $b+ \bar{b} + \mathbf{u}(B + T \bar{B})$ as $\bar{b}+(b+  \mathbf{u}B) + \mathbf{u}T \bar{B}$, 
from the Theorems (\ref{th:iso1}) and (\ref{th:iso2}), we can deduce the 
following theorem. 
\begin{theorem} \label{th:last}
The $\mathsf{E}^0$-term of the spectral sequence is isomorphic to the complex 
\[ \mathsf{E}^0_{pq}=(\mathsf{C}_q(\mathcal{H},\mathsf{C}_p((A^{op})^{\natural}_{\mathcal{H}^{cop}}) \boxtimes \mathsf{W} ),\delta), \]
and the $\mathsf{E}^1$-term is  \[
\mathsf{E}^1_{pq}= (H_q( \mathcal{H} , \mathsf{C}_p((A^{op})^{\natural}_{\mathcal{H}^{cop}}) \boxtimes \mathsf{W} ) , \mathfrak{b} + \mathbf{u} \mathfrak{B}). \]
The $\mathsf{E}^2$-term of the spectral sequence is
\[ \mathsf{E}^2_{pq} = HC_p( H_q(\mathcal{H}, \mathsf{C}_{\bullet}((A^{op})^{\natural}_{\mathcal{H}^{cop}}));\mathsf{W}), \]
the cyclic homology of the cyclic module $H_p(\mathcal{H},\mathsf{C}_{\bullet}((A^{op})^{\natural}_{\mathcal{H}^{cop}}))$ with coefficieints in $\mathsf{W}.$ 
\end{theorem}
If $S$ is invertible, by filtering the complex $\mathsf{C}_{\bullet} (A^{op} \natural \mathcal{H}^{cop}),$ 
by a similar method we conclude
\begin{prop}
The $\mathsf{E}^0$-term of the spectral sequence is isomorphic to the complex 
\[ \mathsf{E}^0_{pq}=(\mathsf{C}_q(\mathcal{H},\mathsf{C}_p(A^{\natural}_{\mathcal{H}}) \boxtimes \mathsf{W} ),\delta), \]
and the $\mathsf{E}^1$-term is  \[
\mathsf{E}^1_{pq}= (H_q( \mathcal{H} , \mathsf{C}_p(A^{\natural}_{\mathcal{H}}) \boxtimes \mathsf{W} ) , 
\mathfrak{b} + \mathbf{u} \mathfrak{B}). \]
The $\mathsf{E}^2$-term of the spectral sequence is
\[ \mathsf{E}^2_{pq} = HC_p( H_q(\mathcal{H}, \mathsf{C}_{\bullet}(A^{\natural}_{\mathcal{H}}));\mathsf{W}), \]
the cyclic homology of the cyclic module $H_p(\mathcal{H},\mathsf{C}_{\bullet}(A^{\natural}_{\mathcal{H}}))$ with coefficieints in $\mathsf{W}.$ 
\end{prop}

We give an application of the above spectral sequence. Let $k$ be a field. A Hopf algebra $\mathcal{H}$ over $k$ is called semisimple
if it is semisimple as an algebra, i.e., every injection $ N \hookrightarrow M$ of $\mathcal{H}$-modules splits in the category
of $\mathcal{H}$-modules. It is known that semisimple Hopf algebras are finite dimensional and $\mathcal{H}$ is
semisimple if and only if there is a right integral $t \in \mathcal{H}$ with $\epsilon(t)=1$~\cite{sw69}. Recall that
a right integral in $\mathcal{H}$ is an element $t \in \mathcal{H}$ such that for all $h \in \mathcal{H}$, $th=\epsilon(h)t.$
Now, it is clear that if $\mathcal{H}$ is semisimple, then the functor of coinvariants $(\mathcal{H}\!\!-\!\!mod) \rightarrow (k\!\!-\!\!mod), M \rightarrow M_{\mathcal{H}}$
is exact. It follows that the higher derived functors are zero. In other words, for any left $\mathcal{H}$-module $M$, we have 
$H_0(\mathcal{H},M)= M_{\mathcal{H}}$ and $H_i(\mathcal{H},M)=0$ for $i > 0.$ Alternatively, using a right integral $t \in \mathcal{H}$ with
$\epsilon(t)=1$, we have the following homotopy operator $h:  \mathcal{H}^{\otimes n}  \otimes M \rightarrow   \mathcal{H}^{\otimes (n+1)} \otimes M, n \ge 0,$
$$ h( h_1 \otimes \dots \otimes h_n \otimes m ) =  t \otimes h_1 \otimes \dots \otimes h_n  \otimes m.$$
One can easily check that $\delta h + h \delta =  id$.
\begin{prop}
Let $\mathcal{H}$ be a semisimple Hopf algebra. Then there is a natural isomorphism of cyclic homology groups
$$HC_{\bullet} (A \rtimes \mathcal{H}; \mathsf{W})= HC_{\bullet} (C^{\mathcal{H}}_{\bullet}(A); \mathsf{W}), $$
where $C^{\mathcal{H}}_{\bullet}(A)= H_0(\mathcal{H}, A^{\mathcal{H}}_{\natural})$ is the cyclic module of equivariant chains.
\end{prop}
\begin{proof}
Since $\mathcal{H}$ is semisimple, we have $\mathsf{E}^1_{pq}=0$ for $q>0$ and the spectral sequence collapses. 
The first row of $\mathsf{E}^1$ is exactly 
$H_{\bullet}(\mathcal{H},A^{\mathcal{H}}_{\natural})=C_{\bullet}^{\mathcal{H}}(A).$
\end{proof}
We give an application to cyclic homology of invariant subalgebras. Let $A^{\mathcal
{H}}=\{ a \in A; ha=\epsilon(h)a \;\; \forall h \in \mathcal{H}\}.$ It is easy to check that $A^{\mathcal{H}}$ is a subalgebra of $A$. 
We define a map
\begin{eqnarray*}
& &\chi: C_{\bullet}(A^{\mathcal{H}}) \rightarrow C_{\bullet}^{\mathcal{H}}(A)\\
& &\chi (a_0,\dots,a_n)=(1 \mid a_0,\dots,a_n).
\end{eqnarray*} 
One can check that $\chi$ is a morphism of cyclic modules. It is shown in~\cite{sMo82}, corollary 4.5.4, that under suitable
conditions $A^{\mathcal{H}}$ and $A \rtimes \mathcal{H}$ are Morita equivalent. Thus, in this case, Theorem~\ref{th:last} gives
an spectral sequence converging to $HC_{\bullet} (A^{\mathcal{H}};\mathsf{W}).$ 

We give another application of Proposition 5.2. Let $A$ be an algebra over a field of characteristic zero and 
$\delta : A \rightarrow A$ a locally nilpotent derivation (i.e. for any $a \in A$, there is an integer 
$n > 0$ with $\delta^n(a)=0$). Let $\mathcal{H}=k[x]$ be the universal enveloping algebra of free Lie algebra 
with one generator $x$. Then $A$ is an $\mathcal{H}$-module algebra and we consider the crossed product 
$A \rtimes \mathcal{H}.$ In~\cite{tm85}, using an infinitesimal analogue of Takesaki's duality theorem 
(due to Y. Nouaz\'{e} and P. Gabriel), it is shown that $ HP^{\ast}(A \rtimes \mathcal{H}) \simeq HP^{\ast}(A).$
We show that this result easily follows from Proposition 5.2. In fact, in this case the
$\mathsf{E}^1$-term has only two non-zero rows given by $\mathsf{E}^1_{p0}=(A^{\natural}_{\mathcal{H}})^{\mathcal{H}}$ and
$\mathsf{E}^1_{p1}=(A^{\natural}_{\mathcal{H}})_{\mathcal{H}}$. We obtain a long exact sequence
\[
\dots \rightarrow HC_p((A^{\natural}_{\mathcal{H}})_{\mathcal{H}}) \rightarrow HC_p(A \rtimes {\mathcal{H}}) \rightarrow 
HC_{p+1}((A^{\natural}_{\mathcal{H}})^{\mathcal{H}})\rightarrow 
HC_{p-1}((A^{\natural}_{\mathcal{H}})_{\mathcal{H}}) 
\rightarrow \dots
\]
(cf.~\cite{cw94}, page 124). Now we dualize to cohomology and take $\underrightarrow{lim}$. Since $\delta$ is nilpotent, in 
$A^{\natural}_{\mathcal{H}}=\bigoplus_{n \ge 0}x^n \otimes A^{\natural}$, the terms involving $x^n, n \ge 1,$ do not contribute to periodic cyclic cohomology 
(Proposition 5.1). We then have a $6$-term exact sequence

\[
\begin{CD}
HP^{odd}(A) @<< < HP^{odd}(A \rtimes \mathcal{H}) @<0 << HP^{ev}(A) \\
@V 0 VV  @.                        @AA 0 A \\
HP^{odd}(A) @>0 >> HP^{ev}(A \rtimes \mathcal{H}) @> >> HP^{ev}(A) 
\end{CD}
\]

The vertical arrows correspond to action of the derivation and are therefore zero. The other zero arrows are 
cup product with the generator of $HP^{odd}(\mathcal{H})=0$. This proves our claim. 

\section*{Acknowledgments}

We would like to thank C. Voigt and B. Rangipour for carefully reading an earlier version of this paper and making
useful remarks. We would also like 
to heartilly thank the referee for noticing that our original  formulas in 
Proposition 5.1  
always simplify, irrespective of $\mathcal{H}$ being cocommutative or not, provided we use 
a more appropriate action~(\ref{eq:ac1}) suggested by him. The coinvariant spaces do not 
change.
%********************************** Bibliography ****************************

\end{document}